\documentclass{amsart}
\usepackage{placeins,float}
\usepackage{amsmath,amssymb,amsthm,mathtools,booktabs,array,microtype}
\usepackage[svgnames]{xcolor}
\usepackage{hyperref}
\hypersetup{
  colorlinks=true,
  allcolors=Navy,
  breaklinks=true,
  pdfauthor={Zhichao Tang and Xiang Fan},
  pdftitle={Wendt's Determinant and Finiteness of Prime Characteristics for Permutation Binomials},
  pdfkeywords={permutation binomials, finite fields, prime characteristics, Wendt's determinant, inverse power sums}
}
\newtheorem{theorem}{Theorem}[section]
\newtheorem{lemma}[theorem]{Lemma}
\newtheorem{proposition}[theorem]{Proposition}
\newtheorem{corollary}[theorem]{Corollary}
\theoremstyle{remark}
\newtheorem{remark}[theorem]{Remark}
\title[Wendt's Determinant and Permutation Binomials]{Wendt's Determinant and Finiteness of Prime Characteristics for Permutation Binomials}
\author{Zhichao Tang}
\author{Xiang Fan}
\address{School of Mathematics, Sun Yat-sen University, Guangzhou 510275, China}
\email[Xiang Fan]{fanx8@mail.sysu.edu.cn}
\date{}
\begin{document}
\begin{abstract}
For $f(X)=X^r(X^s+a)\in\mathbb F_{p^m}[X]$, with
$m,r,s\geqslant1$ and $a\in\mathbb F_{p^m}^\times$, put
$g=\gcd(s,p^m-1)$. We prove that, for each fixed $g$,
only finitely many prime characteristics $p$ admit such
permutation binomials, even as $m,r,s,a$ vary.
If $6\nmid g$, every such $p$ divides Wendt's determinant
$W_g=\operatorname{Res}_X(X^g-1,(X+1)^g-1)$, which is nonzero.
When $6\mid g$, this determinant vanishes; inverse power sums
and iterated resultants yield an explicit nonzero integer
divisible by every such $p>2g+1$.
We determine the exact admissible characteristics for
$1\leqslant g\leqslant10$, obtain a further divisibility
condition over prime fields, and show that in prime fields
$g=11$ occurs precisely when $p=23$.
\end{abstract}
\maketitle

\section{Introduction}
Let $p$ be prime, $m\geqslant1$, and $q:=p^m$. Consider
\begin{equation}\label{eq:f}
 f(X):=X^r(X^s+a),\qquad r,s\geqslant1,\quad
 a\in\mathbb{F}_{p^m}^\times,\qquad g:=\gcd(s,p^m-1).
\end{equation}
We fix $g$ and study the characteristics $p$ for which $f$
permutes $\mathbb F_{p^m}$, allowing $m,r,s,a$ to vary.
For every integer $g\geqslant1$, define
\[
 \mathcal{P}_g:=\left\{p\text{ prime}:\ \begin{aligned}
 &\exists\text{ integers }m,r,s\geqslant1,\quad
   a\in\mathbb{F}_{p^m}^\times\text{ such that}\\
 &X^r(X^s+a)\text{ permutes }\mathbb{F}_{p^m},\quad
  \gcd(s,p^m-1)=g
 \end{aligned}\right\}.
\]
There is no bound on the extension degree $m$.
We prove that $\mathcal P_g$ is finite for every fixed $g$.
The proof also gives explicit obstructions when $6\mid g$,
in which case Wendt's determinant (also known as
Wendt's binomial circulant determinant) vanishes.

For every $g\geqslant1$, define \emph{Wendt's determinant} by
\[
 W_g:=\operatorname{Res}_X(X^g-1,(X+1)^g-1).
\]
As usual, $\operatorname{Res}_X(F,G)$ denotes the resultant
of $F,G\in\mathbb Z[X]$.
When $6\mid g$, also define
\[
 E_g:=\operatorname{Res}_X\left(\frac{X^g-1}{X^2-X+1},(X-1)^g-1\right).
\]
For $g\geqslant2$, define $T_{g,k}\in\mathbb Z[X,Y]$ recursively by
\begin{equation}\label{eq:Trec}
 T_{g,1}:=1,\qquad
 T_{g,k+1}:=\bigl[(X-gk)(Y-1)+kY\bigr]T_{g,k}
             -Y(Y-1)\partial_YT_{g,k}.
\end{equation}
For even $g\geqslant2$, set
\[
 R_g:=\operatorname{Res}_Y\!\left(Y^2-Y+1,
 \operatorname{Res}_X\bigl(T_{g,g}(X,Y),
                   T_{g,2g}(2X+1,Y)\bigr)\right)\in\mathbb Z.
\]
The inner resultant lies in $\mathbb Z[Y]$; in particular,
$R_g\in\mathbb Z$.

\begin{theorem}[Characteristic finiteness]\label{thm:finite}
For every integer $g\geqslant1$, the set $\mathcal{P}_g$ is finite.
If $6\nmid g$, Wendt's determinant $W_g$ is nonzero,
and every $p\in\mathcal{P}_g$ divides $W_g$.
If $6\mid g$, the integers
$E_g$ and $R_g$ are nonzero and satisfy
\begin{equation}\label{eq:maincontain}
 \mathcal{P}_g\subseteq\{p:p\leqslant2g+1\}
 \cup\{p:p\mid E_gR_g\}.
\end{equation}
Moreover, $R_g$ is odd.
\end{theorem}

\begin{theorem}[Possible characteristics for small $g$]\label{thm:spectrum}
For $1\leqslant g\leqslant10$, the sets $\mathcal P_g$ are
\[
 \begin{array}{@{}cl@{}}
\toprule
 g&\mathcal P_g\\\midrule
 1&\varnothing\\
 2&\{3\}\\
 3&\{2,7\}\\
 4&\{3,5\}\\
 5&\{11\}
 \\
\bottomrule
\end{array}
 \qquad\qquad
 \begin{array}{@{}cl@{}}
\toprule
 g&\mathcal P_g\\\midrule
 6&\{5,7,13,19,31\}\\
 7&\{2,29\}\\
 8&\{3,5,17\}\\
 9&\{2,7,19,37\}\\
 10&\{3,11,31\}
 \\
\bottomrule
\end{array}
\]
\end{theorem}

For prime fields, Masuda and Zieve~\cite[Theorem~1.1]{MZ07}
excluded $g\in\{2,4\}$ when $p>5$. These exclusions also
follow from $p\mid W_g$, since
$W_2=-3$ and $W_4=-3\cdot5^3$. Theorem~\ref{thm:spectrum}
also determines the admissible characteristics for
$g\in\{2,4\}$ over extensions of arbitrary degree. Masuda and Zieve
also proved $p-1\leqslant g(g+1)$~\cite[Theorem~1.3]{MZ};
Ayad, Belghaba, and Kihel improved this to
$p-1\leqslant g(g-1)$~\cite[Theorem~1.7]{ABK14}.
For earlier prime-field results, see Wan~\cite{Wan87}.

Over extension fields, these inequalities for $p-1$
do not extend to $q-1$:
$X^{10}+3X$ permutes $\mathbb F_{7^3}$~\cite{MZ},
where $g=9$ and $q-1=342>g(g+1)$.
Ayad, Belghaba, and Kihel~\cite[Corollary~3.6]{ABK14}
proved a conditional descent result: if $f$ in~\eqref{eq:f}
permutes $\mathbb F_{p^m}$ and
$\gcd(\ell(\ell-1),m)=1$ for every prime $\ell\mid g$,
then either $g=p-1$ or there is an associated
permutation binomial over $\mathbb F_p$.
In particular, their Corollary~3.7(v) gives
$p\in\{7,13,19,31\}$ for $g=6$ when $\gcd(m,6)=1$.
Theorem~\ref{thm:spectrum} removes the restriction on $m$ and
also includes $p=5$, realized over $\mathbb F_{25}$.
For $g>1$, their later necessary condition
$p\not\equiv1\pmod{g^2}$~\cite[Theorem~2.3]{ABK15}
excludes a residue class but does not imply finiteness.
This condition is therefore complementary to finiteness.
Akbary and Wang~\cite{AkbaryWang} instead used the parameter
$\ell=(q-1)/\gcd(v-r,q-1)$ for binomials $X^r+X^v$;
fixing $\ell$ is different from fixing $g$ as $q$ varies.

The product of the nonzero values of~\eqref{eq:f}
gives a divisibility condition involving Wendt's determinant
\cite{Wendt,FeeGranville}. For related circulant methods,
see~\cite[\S3]{ParkLee}.

The determinant vanishes precisely when $6\mid g$.
Fee and Granville~\cite[\S1]{FeeGranville} studied its
cyclotomic modification $E_g$.
For $6\mid g$, we obtain a second obstruction from two inverse
power sums. Eliminating their common parameter, then taking a
resultant with $Y^2-Y+1$, gives $R_g$. A reduction modulo $2$
proves that $R_g\ne0$, leaving only finitely many possible
characteristics.
Over prime fields, inverse power sums yield an additional
divisibility condition; see Theorem~\ref{thm:primefield-content}.
The prime-field case $g=11$ is settled by summing the values
of the quotient map (Corollary~\ref{cor:primefield-g11}).

The classification for $g\leqslant10$ uses exact finite
calculations; some prime-field cases were known already
\cite[Corollary~2.5]{MZ}. Theorem~\ref{thm:spectrum}
classifies characteristics rather than individual binomials.
For seven parameter pairs, Corollary~\ref{cor:fan-seven}
gives a permutation criterion using~\cite[Corollary~6.1]{Fan26}.
The other characteristic pairs are not covered by that
criterion. Small-order factorizations of Wendt's determinant are tabulated
in~\cite{MichonWendt} and, for even orders, in~\cite{FeeGranville}.
The remaining finite certificates are recorded in
Appendices~\ref{app:exclusions}--\ref{app:witnesses} and
checked by the accompanying script.

\section{The product and quotient obstructions}
\begin{lemma}\label{lem:product}
Let $\mathbb F_q$ be a finite field, and let $r,s\geqslant1$
and $a\in\mathbb F_q^\times$. Suppose
$f(X)=X^r(X^s+a)$ permutes $\mathbb F_q$.
Put $g=\gcd(s,q-1)$, $t=(q-1)/g$, and $\zeta=(-a)^t$.
Then $g>1$, $\gcd(r,g)=1$, and
\begin{equation}\label{eq:prod}
 \zeta^g=1,\qquad \zeta\ne1,\qquad
 (\zeta-1)^g=(-1)^{r-1}.
\end{equation}
\end{lemma}
\begin{proof}
If $g=1$, the map $x\mapsto x^s$ is bijective on
$\mathbb F_q^\times$. Hence $x^s=-a$ has a nonzero solution,
contrary to the injectivity of $f$, since $f(0)=0$.

For $x\ne0$ and $\xi^s=1$, we have
$f(\xi x)=\xi^r f(x)$. The kernel of $x\mapsto x^s$ on $\mathbb F_q^\times$ has
order $g$. As $f(x)\ne0$, injectivity on each fiber forces
$\xi\mapsto\xi^r$ to be injective on this kernel. Thus
$\gcd(r,g)=1$.

The image of $x\mapsto x^s$ consists of the $t$-th
roots of unity, each attained $g$ times. If $(-a)^t=1$, then $x^s=-a$
for some $x\ne0$, contrary to $f(x)\ne f(0)$.
Thus $\zeta\ne1$, and
$\zeta^g=(-a)^{q-1}=1$.
Since $f$ permutes $\mathbb F_q$ and $f(0)=0$,
$\prod_{x\ne0}f(x)=\prod_{x\ne0}x=-1$.
Dividing by $\prod_{x\ne0}x^r=(-1)^r$ gives
\[
 (-1)^{1-r}=\prod_{x\ne0}(x^s+a)
 =\left(\prod_{y^t=1}(y+a)\right)^g
 =\bigl((-1)^t(\zeta-1)\bigr)^g.
\]
The last equality follows by evaluating $Y^t-1$ at
$Y=-a$. Since $gt=q-1$, we have $(-1)^{gt}=1$ in
$\mathbb F_q$. This proves~\eqref{eq:prod}.
\end{proof}

\begin{proposition}[Wendt obstruction]\label{prop:wendt}
Let $\mathbb F_q$ be a finite field of characteristic $p$.
If $f(X)=X^r(X^s+a)$ permutes $\mathbb F_q$, where
$r,s\geqslant1$ and $a\in\mathbb F_q^\times$, then
$p\mid W_g$ and $p\nmid g$, with $g=\gcd(s,q-1)$.
\end{proposition}
\begin{proof}
Since $g\mid q-1$, we have $p\nmid g$.
Put $t=(q-1)/g$ and $\zeta=(-a)^t$ as in Lemma~\ref{lem:product}.
Set $\eta=\zeta$ if $g$ is even or $r$ is odd, and
$\eta=\zeta^{-1}$ otherwise. In the latter case $g$ is
odd and $r$ is even, so~\eqref{eq:prod} gives
$(\eta-1)^g=(-1)^g(\zeta-1)^g=1$.
In the former case this follows directly from~\eqref{eq:prod}.
In both cases $\eta^g=(\eta-1)^g=1$.
The polynomials $X^g-1$ and $(X+1)^g-1$ therefore
have the common root $\eta-1$ in $\mathbb F_q$.
Their resultant $W_g$ vanishes in $\mathbb F_q$, hence
$p\mid W_g$.
\end{proof}

The vanishing criterion for Wendt's determinant is classical
\cite{Lehmer35}; we include a short proof.
\begin{lemma}\label{lem:wendt-zero}
For every integer $g\geqslant1$, one has $W_g=0$
if and only if $6\mid g$.
\end{lemma}
\begin{proof}
Suppose $W_g=0$. Then $X^g-1$ and $(X+1)^g-1$
have a common root $z\in\mathbb C$. Thus
$|z|=|z+1|=1$, so $\operatorname{Re}z=-1/2$.
Hence $z$ is a primitive cube root of unity and $z+1$
a primitive sixth root. In particular, $6\mid g$.
Conversely, if $6\mid g$, either primitive cube root
$z$ satisfies $z^g=(z+1)^g=1$. Hence $W_g=0$.
\end{proof}
\begin{remark}[Bounds for prime divisors]\label{rem:wendt-bounds}
Suppose $g\geqslant2$ and $6\nmid g$. Then every prime $p\mid W_g$
satisfies
\begin{equation}\label{eq:wendt-prime-bound}
 p\mid 2^g-1\quad\text{or}\quad p\leqslant3^{g/2};
 \qquad\text{in particular,}\quad p\leqslant2^g-1.
\end{equation}
Indeed, for a primitive $g$-th root $\omega$, the integer products
$N_{a,b}=\prod_{j=0}^{g-1}(1+\omega^{aj}-\omega^{bj})$
account for all prime divisors of
$W_g=\prod_{a,b\bmod g}(1+\omega^a-\omega^b)$:
reduce this identity modulo a prime ideal above $p$.
If none of $a,b,a-b$ is zero modulo $g$, orthogonality
and the arithmetic--geometric mean inequality give
$|N_{a,b}|\leqslant3^{g/2}$; otherwise $N_{a,b}$ is a unit
or a power of $2^d-1$, where $d\mid g$.

For even $g$ with $3\nmid g$, Simalarides
\cite[Theorem~3]{Simalarides02} obtained the sharper assertion
that,
with $\rho=2.2487338$ and $L_n$ the Lucas numbers
($L_0=2$, $L_1=1$, $L_{n+2}=L_{n+1}+L_n$),
we have
\[
 p>\rho^{g/4}\quad\Longrightarrow\quad
 \begin{cases}
 p\mid (2^g-1)L_{g/2}, & g\equiv2\pmod4,\\
 p\mid 2^{g/2}+1, & g\equiv0\pmod8,\\
 (g,p)\in\{(4,3),(4,5),(20,61)\}, & g\equiv4\pmod8.
 \end{cases}
\]
Consequently $p\leqslant2^{g/2}+1$ when $g\equiv0\pmod8$,
and $p\leqslant\rho^{g/4}$ when $g\equiv4\pmod8$
and $g\geqslant28$.

For odd $g$ with $3\nmid g$, multiplicativity of resultants
gives $W_g\mid W_{2g}$, since $X^g-1$ divides
$X^{2g}-1$ and $(X+1)^g-1$ divides $(X+1)^{2g}-1$.
Applying Simalarides~\cite[Theorem~3]{Simalarides02}
to the even order $2g\equiv2\pmod4$ therefore yields
\[
 p\mid (2^{2g}-1)L_g\quad\text{or}\quad
 p\leqslant\rho^{g/2}.
\]
When $g$ is odd and $3\mid g$, this argument does not apply
because $W_{2g}=0$; the general bound~\eqref{eq:wendt-prime-bound}
still holds.
\end{remark}

The reduction of a permutation binomial to a map on a
root-of-unity subgroup appears in
\cite[Theorem~1.2]{WanLidl} and \cite[Lemma~2.1]{Zieve};
see also the quadratic-extension case in
\cite[Theorem~1.12]{HouPL}.
We use the following precise form.
For each positive divisor $d$ of $q-1$, write
$\mu_d:=\{x\in\mathbb F_q^\times:x^d=1\}$.

\begin{lemma}[Quotient normalization]\label{lem:quotient}
Let $\mathbb F_q$ be a finite field, $r,s\geqslant1$,
$a\in\mathbb F_q^\times$, and $f(X)=X^r(X^s+a)$.
Put $g=\gcd(s,q-1)$, $t=(q-1)/g$, and $u=s/g$,
so that $\gcd(u,t)=1$.
Choose $v\in\mathbb Z$ with $uv\equiv1\pmod t$, and let
$0\leqslant b<t$ satisfy $b\equiv rv\pmod t$.
Then $f$ permutes $\mathbb F_q$ if and only if $\gcd(r,g)=1$ and
$\psi(X):=X^b(X+a)^g$ permutes $\mu_t$.
\end{lemma}
\begin{proof}
The map $x\mapsto x^g$ induces a bijection
$\mathbb F_q^\times/\mu_g\cong\mu_t$. For $y=x^g\in \mu_t$,
\[
 f(x)^g=y^r(y^u+a)^g=\psi(y^u),
\]
since $ub\equiv r\pmod t$; moreover, $y\mapsto y^u$ permutes
$\mu_t$. For $\xi\in\mu_g$ we have $f(\xi x)=\xi^r f(x)$.
If $f$ permutes $\mathbb F_q$, then $\psi$ permutes
$\mu_t$, and injectivity on each $\mu_g$-coset forces
$\gcd(r,g)=1$. Conversely, these two conditions ensure
that $f$ permutes the nonzero $\mu_g$-cosets and is
bijective on each coset. Together with $f(0)=0$, this
shows that $f$ permutes $\mathbb F_q$.
\end{proof}

For $z\in\mathbb F_q^\times$, write
$\psi_{b,z}(Y):=Y^b(Y-z)^g$ on $\mu_t$.
If $b'\equiv-b-g\pmod t$, then
\begin{equation}\label{eq:quotient-inversion}
 \psi_{b,z}(y^{-1})=(-z)^g\psi_{b',z^{-1}}(y)
 \qquad(y\in \mu_t).
\end{equation}
Since $((-z)^g)^t=1$, inversion on $\mu_t$ and
multiplication by $(-z)^g$ are both permutations.
Thus the two maps permute $\mu_t$ simultaneously.

\section{Sixth roots and inverse power sums}\label{sec:moments}
Throughout this section, let $\mathbb F_q$ be a finite field
of characteristic $p$, and let
$f(X)=X^r(X^s+a)$ with $r,s\geqslant1$ and
$a\in\mathbb F_q^\times$. Set $g=\gcd(s,q-1)$,
$t=(q-1)/g$, and $u=s/g$, so that $\gcd(u,t)=1$.
Choose $v\in\mathbb Z$ with $uv\equiv1\pmod t$,
and $0\leqslant b<t$ with $b\equiv rv\pmod t$.
Write $\mu_t=\{y\in\mathbb F_q^\times:y^t=1\}$,
$\psi(X)=X^b(X+a)^g$, $\zeta=(-a)^t$,
and $\tau=-g^{-1}\in\mathbb F_p$ (well defined since $g\mid q-1$).

Let $\Phi_6(X):=X^2-X+1$ be the sixth cyclotomic polynomial.
If $\psi$ permutes $\mu_t$, the inverse power sums
$\sum_{y\in\mu_t}\psi(y)^{-j}$ vanish for $1\leqslant j<t$.
For $6\mid g$, recall that
\[
 E_g=\operatorname{Res}_X\left(
 \frac{X^g-1}{\Phi_6(X)},(X-1)^g-1\right).
\]
The complex roots of $\Phi_6$ have multiplicative order $6$.
For $6\mid g$, $E_g$ is the cyclotomic modification of
Wendt's determinant used by Fee and Granville~\cite[\S1]{FeeGranville};
the unmodified $W_g$ vanishes.
\begin{lemma}\label{lem:residual}
We have $E_g\ne0$. If $f$ permutes $\mathbb F_q$ and
$p\nmid E_g$, then $\Phi_6(\zeta)=0$.
\end{lemma}
\begin{proof}
The only common complex roots of $X^g-1$ and
$(X-1)^g-1$ are the primitive sixth roots. Removing
$\Phi_6$ therefore leaves a nonzero resultant.
As $g$ is even, $\gcd(r,g)=1$ forces $r$ odd;
by~\eqref{eq:prod}, $(\zeta-1)^g=1$.
If $\Phi_6(\zeta)\ne0$, the polynomials
$(X^g-1)/\Phi_6(X)$ and $(X-1)^g-1$ have the common
root $\zeta$ in $\mathbb F_q$. Their resultant $E_g$
vanishes in $\mathbb F_q$, hence $p\mid E_g$.
\end{proof}

The following inverse-power-sum identities hold without the
assumption $6\mid g$.
For integers $1\leqslant j<t$, put
\[
 b_j\equiv jb\pmod t,\quad 0\leqslant b_j<t,\qquad
 U_j:=t-b_j\in\{1,\ldots,t\}.
\]
For indeterminates $T,X,Y$, define
$P_k(T,X,Y)\in\mathbb Z[T,X,Y]$ recursively by
\begin{equation}\label{eq:prec}
 P_1(T,X,Y):=1,\qquad
 P_{k+1}:=\bigl[(X-k)(Y-1)-kTY\bigr]P_k
       +TY(Y-1)\partial_YP_k.
\end{equation}

For a formal power series $F(w)=\sum_{j\geqslant0}c_jw^j$
with coefficients in $\mathbb Q(T,X,Y)$, write
$[w^j]F(w):=c_j$, and use the binomial expansion
$(1+w)^A=\sum_{h\geqslant0}\binom{A}{h}w^h$.

\begin{lemma}\label{lem:coeff}
For every integer $k\geqslant1$, the identity
\begin{equation}\label{eq:hasse}
 P_k(T,X,Y)=(k-1)!(Y-1)^k[w^{k-1}]
 \frac{(1+w)^{X-1}}{Y(1+w)^T-1}
\end{equation}
holds in $\mathbb Q(T,X,Y)$.
\end{lemma}
\begin{proof}
Put $A(w)=(1+w)^{X-1}/(Y(1+w)^T-1)$.
Since
$(1+w)\partial_w A=(X-1)A+TY\partial_Y A$,
comparing coefficients shows that
$(k-1)!(Y-1)^k[w^{k-1}]A$ satisfies
\eqref{eq:prec}, with initial value $1$.
\end{proof}

\begin{lemma}[Inverse-power-sum condition]\label{lem:moment}
Assume that $\psi$ permutes $\mu_t$. Then
\begin{equation}\label{eq:momroot}
 P_{gj}(\tau,U_j,\zeta)=0
 \qquad\text{for all }1\leqslant j<t.
\end{equation}
\end{lemma}
\begin{proof}
Fix $1\leqslant j<t$. Since $\psi$ permutes $\mu_t$,
\[
 0=\sum_{y\in\mu_t}\psi(y)^{-j}
  =\sum_{y^t=1}\frac{y^{U_j}}{(y+a)^{gj}}.
\]
Suppose first that $gj\leqslant p$. For $1\leqslant U\leqslant t$,
the rational-function identity
\begin{equation}\label{eq:partialfrac}
 \sum_{y^t=1}\frac{y^U}{Z-y}
 =\frac{t Z^{U-1}}{Z^t-1}
\end{equation}
holds with an indeterminate $Z$. Since $\psi$ permutes
$\mu_t$, we have $(-a)^t\ne1$; otherwise $-a\in\mu_t$
and $\psi(-a)=0$. Also $p\nmid t$ since $gt=q-1$.
By Lemma~\ref{lem:coeff}, specializing $T=t$, $X=U$,
$Y=Z^t$ in the Taylor expansion of the right side
of~\eqref{eq:partialfrac} gives, for $1\leqslant k\leqslant p$,
\[
 \left(\frac d{dZ}\right)^{k-1}\frac{tZ^{U-1}}{Z^t-1}
 =\frac{tZ^{U-k}}{(Z^t-1)^k}P_k(t,U,Z^t).
\]
Differentiating the left side $k-1$ times gives
$-(k-1)!\sum_{y^t=1}y^U/(y-Z)^k$.
Taking $U=U_j$, $k=gj$ and $Z=-a$, this derivative
vanishes by the inverse-power-sum identity.
All other factors are nonzero at $Z=-a$, whence
$P_{gj}(t,U_j,(-a)^t)=0$.

If $gj>p$, Lemma~\ref{lem:coeff} and a geometric-series
expansion in $w$ show that
$(Y-1)^{gj}[w^{gj-1}](1+w)^{U_j-1}/(Y(1+w)^t-1)$
belongs to $\mathbb Z[Y]$, since $t,U_j$ are positive integers.
Consequently $P_{gj}(t,U_j,Y)\in(gj-1)!\mathbb Z[Y]$,
so this polynomial vanishes modulo $p$. Since $t\equiv-g^{-1}=\tau\pmod p$ and
$\zeta=(-a)^t$, both cases give~\eqref{eq:momroot}.
\end{proof}

\section{Nonvanishing and characteristic finiteness}\label{sec:odd}
Recall the polynomials $T_{g,k}$ from~\eqref{eq:Trec}.
Before specializing $Y$, the recurrence~\eqref{eq:prec}
gives
\begin{equation}\label{eq:Tscale}
 T_{g,k}(X,Y)=g^{k-1}P_k(-1/g,X/g,Y).
\end{equation}
For even $g\geqslant2$ and $c\in\{0,1\}$, write
\[
 H_{g,c}(Y):=\operatorname{Res}_X\bigl(T_{g,g}(X,Y),
 T_{g,2g}(2X+1-c,Y)\bigr)\in\mathbb Z[Y].
\]
Thus $R_g=\operatorname{Res}_Y(\Phi_6,H_{g,0})$.

\begin{lemma}[Reflection symmetry]\label{lem:reflection}
For $k\geqslant1$,
\begin{equation}\label{eq:preflect}
 P_k(T,X,Y)=Y^{k-1}P_k(T,k+T-X,Y^{-1}).
\end{equation}
For every even $g\geqslant2$,
\[
 \operatorname{Res}_Y(\Phi_6,H_{g,1})=R_g.
\]
\end{lemma}
\begin{proof}
Put $n=k-1$. The coefficient-change identity
\[
 [w^n]A(w)=(-1)^n[v^n](1+v)^{n-1}
 A\!\left(-\frac{v}{1+v}\right)
\]
is verified on monomials $A(w)=w^j$. Applied to the
coefficient in~\eqref{eq:hasse}, it gives
$(-1)^{n+1}[v^n](1+v)^{n+T-X}/((1+v)^T-Y)$.
Equation~\eqref{eq:hasse} now proves~\eqref{eq:preflect}. By~\eqref{eq:Tscale},
\[
 T_{g,k}(X,Y)=Y^{k-1}T_{g,k}(gk-1-X,Y^{-1}).
\]
Let $d_1=g-1$ and $d_2=2g-1$, both odd. Replacing
$X$ by $g^2-1-X$ in the two inputs of $H_{g,c}$ gives
\[
 H_{g,c}(Y)=(-1)^{d_1d_2}Y^{2d_1d_2}H_{g,1-c}(Y^{-1}).
\]
The two roots of $\Phi_6$ are reciprocal and have
product $1$. Evaluation at these roots proves the claimed
identity of integer resultants.
\end{proof}

\begin{theorem}[Oddness]\label{thm:odd}
For every even $g\geqslant2$, the integer $R_g$ is positive and odd.
\end{theorem}
\begin{proof}
Reduce the full recurrence~\eqref{eq:Trec} modulo $2$
before setting $X=1$. Starting from
$T_{g,2}(1,Y)\equiv1$, induction gives
\[
 T_{g,k}(1,Y)\equiv(Y+1)^{k-2}\pmod2
 \qquad(k\geqslant2).
\]
Let $z\in\mathbb F_4$ satisfy $z^2+z+1=0$.
The first input of $H_{g,0}$ has degree $d_1=g-1$ in $X$
and leading coefficient $(z-1)^{d_1}\ne0$. The second input,
whose degree before reduction is $d_2=2g-1$, reduces
modulo $2$ to the nonzero constant $(z+1)^{2g-2}$.
We specialize the resultant with the degrees $d_1,d_2$
fixed, viewing it as a polynomial in the coefficients. If the second input
reduces to a nonzero constant $b$ and the first has leading
coefficient $a$, the specialization is $a^{d_2}b^{d_1}$.
Indeed, apply the root-product formula to
$b+\varepsilon X^{d_2}$ and then set $\varepsilon=0$. Hence
\[
 H_{g,0}(z)=(z-1)^{d_1d_2}(z+1)^{(2g-2)d_1}\ne0.
\]
As $H_{g,0}$ has integer coefficients, its value at
$z^2$ is $H_{g,0}(z)^2$. Every nonzero element of
$\mathbb F_4$ has cube $1$, and hence
\[
 R_g\equiv H_{g,0}(z)H_{g,0}(z^2)=1\pmod2.
\]
Write $H_{g,0}(Y)\equiv A+BY\pmod{\Phi_6(Y)}$
with $A,B\in\mathbb Z$. By the definition of the resultant,
$R_g=A^2+AB+B^2$. This is positive, since $R_g$ is odd.
\end{proof}

\begin{proof}[Proof of Theorem~\ref{thm:finite}]
Lemma~\ref{lem:product} rules out $g=1$.
For $6\nmid g$, Proposition~\ref{prop:wendt}
and Lemma~\ref{lem:wendt-zero} give the assertion.
Suppose $6\mid g$, $p>2g+1$, and $p\nmid E_g$.
Then $t=(q-1)/g>2$ and $2g<p$.
Lemma~\ref{lem:residual} gives $\Phi_6(\zeta)=0$.
Apply Lemma~\ref{lem:moment} with $j=1,2$ and put
$X_0=gU_1\in\mathbb F_p$ and
$c=\lfloor2b/t\rfloor\in\{0,1\}$.
Since $U_2=2U_1+(c-1)t$ and $gt\equiv-1\pmod p$,
\[
 T_{g,g}(X_0,\zeta)=0,\qquad
 T_{g,2g}(2X_0+1-c,\zeta)=0.
\]
Because $H_{g,c}$ and $\Phi_6$ have the common root $\zeta$,
we have
$p\mid\operatorname{Res}_Y(\Phi_6,H_{g,c})$.
By Lemma~\ref{lem:reflection} this integer is $R_g$,
which is nonzero by Theorem~\ref{thm:odd}.
Thus~\eqref{eq:maincontain} holds.
\end{proof}

\section{Prime-field refinements}\label{sec:primefield}
Over prime fields, the first inverse power sum gives an
additional divisibility condition. We use the recurrence
\eqref{eq:Trec} for all $g\geqslant2$; the identities
\eqref{eq:Tscale} and~\eqref{eq:preflect} remain valid.
For $1\leqslant j<g$, put
\[
 \alpha_{g,j}:=\begin{cases}
gj,&g\text{ even or }j\text{ odd},\\
 g^2-1-gj,&g\text{ odd and }j\text{ even},
 \end{cases}
 \qquad S_g^{\mathrm{pf}}(X):=\prod_{j=1}^{g-1}T_{g,g}(\alpha_{g,j},X).
\]
Put $H_g(X):=(X^g-1)/(X-1)$ and $B_g(X):=(X-1)^g-1$.
Interchanging and translating the factors in the definition
of $W_g$, then separating the factor $X-1$, gives
\[
 \operatorname{Res}_X(H_g(X),B_g(X))=(-1)^{g+1}W_g.
\]
For $Q\in\mathbb Z[X]$, let
$\operatorname{cont}_X(Q)$ be its \emph{content}, the
nonnegative gcd of its integer coefficients. Define
\[
 D_g:=\operatorname{cont}_X
 \operatorname{Res}_Y\bigl(H_g(Y),B_g(Y)+XS_g^{\mathrm{pf}}(Y)\bigr),
\]
which is a nonnegative integer.
\begin{theorem}[A divisibility condition over prime fields]\label{thm:primefield-content}
The integer $D_g$ is nonzero for every $g\geqslant2$.
If $6\nmid g$, it divides $|W_g|$; if $g$ is even, it is odd.
If $X^r(X^s+a)$ permutes $\mathbb F_p$ with $\gcd(s,p-1)=g$,
then either $p=g+1$ or $p\mid D_g$.
\end{theorem}
\begin{proof}
If $t=(p-1)/g=1$, then $p=g+1$.
Now suppose $t\geqslant2$. The prime-field bound
\cite[Theorem~1.7]{ABK14} gives $t\leqslant g-1$.
Choose $0\leqslant b<t$ as in Lemma~\ref{lem:quotient}, and
put $U=t-b$ and $\zeta=(-a)^t$.
Since $1<t$ and $g<p$, Lemma~\ref{lem:moment} applies with $j=1$:
$T_{g,g}(gU,\zeta)=0$ with $1\leqslant U<g$.
For even $g$, the product identity gives $B_g(\zeta)=0$.
For odd $g$, the integer $t$ is even; as $\gcd(u,t)=1$,
both $u=s/g$ and its inverse modulo $t$ are odd.
Thus $b\equiv r\pmod2$, so $U\equiv r\pmod2$.
For odd $U$ take $\zeta$; for even $U$ use
$\zeta^{-1}$ and~\eqref{eq:preflect}, which replaces
$gU$ by $g^2-1-gU$.
Thus $H_g,B_g,S_g^{\mathrm{pf}}$ have a common root
modulo $p$ in either case.

Since $H_g$ is monic, reduction of the resultant modulo any
prime $\ell$ is the product
\[
 \prod_{H_g(\eta)=0}
 \bigl(B_g(\eta)+XS_g^{\mathrm{pf}}(\eta)\bigr)
\]
over a field containing all roots of $H_g$ modulo $\ell$,
with multiplicities; all three polynomials are reduced modulo $\ell$.
The polynomial ring over the splitting field is an integral
domain, so the product vanishes identically in $X$ if and
only if some factor has both coefficients zero. Thus
\[
 \ell\mid D_g
 \iff \deg\gcd_{\mathbb F_\ell[Y]}
 (H_g,B_g,S_g^{\mathrm{pf}})>0.
\]
The common root obtained above therefore gives
$p\mid D_g$.
For $6\nmid g$, the constant coefficient is
$(-1)^{g+1}W_g\ne0$. Hence $D_g\ne0$ and $D_g\mid|W_g|$.
For even $g$, the recurrence at $X=0$ modulo $2$ reads
\[
 T_{g,k+1}(0,Y)=kY T_{g,k}(0,Y)
    +Y(Y+1)\partial_Y T_{g,k}(0,Y).
\]
Starting from $T_{g,2}(0,Y)=Y$, induction yields
$T_{g,g}(gj,Y)\equiv Y(Y+1)^{g-2}\pmod 2$ for every $j$.
Since $H_g(0)=1$ and $B_g(1)\equiv1\pmod2$, the three
polynomials have no common root in characteristic $2$.
Hence $2\nmid D_g$, so $D_g$ is nonzero and odd.
\end{proof}

The content condition may exclude primes that divide
Wendt's determinant. For example, $73\mid W_9$ but
$73\nmid D_9$: over
$\mathbb F_{73}$, we have $\gcd(H_9,B_9)=Y-2$ and
$S_9^{\mathrm{pf}}(2)\equiv2\cdot18\cdot26^2\cdot52\cdot13\cdot36\equiv-1$.
The following corollary is independent of
Theorem~\ref{thm:primefield-content}.

\begin{corollary}[The prime-field case $g=11$]\label{cor:primefield-g11}
Let $p$ be prime. There exist integers $r,s\geqslant1$ and
$a\in\mathbb F_p^\times$ such that
$\gcd(s,p-1)=11$ and $X^r(X^s+a)$ permutes
$\mathbb F_p$ if and only if $p=23$.
\end{corollary}
\begin{proof}
By~\cite[Theorem~1.7]{ABK14}, $p-1\leqslant110$. As
$11\mid p-1$ and $p$ is odd, we have $p\in\{23,67,89\}$.
Put $t=(p-1)/11$, $z=-a$, and take $b$ as in
Lemma~\ref{lem:quotient}. Then $\psi_{b,z}$ permutes $\mu_t$. As $t>1$,
summing over $\mu_t$ gives
\[
 0=\frac1t\sum_{y\in \mu_t}\psi_{b,z}(y)
 =\sum_{\substack{0\leqslant k\leqslant11\\k+b\equiv0\pmod t}}
 \binom{11}{k}(-z)^{11-k}.
\]
Put $\zeta=z^t$, so $\zeta^{11}=1$.
Since $t$ is even, the product of the elements of $\mu_t$
is $-1$. Taking products of $\psi_{b,z}$ yields
$(-1)^b(\zeta-1)^{11}=-1$, or
$(\zeta-1)^{11}=(-1)^{b-1}$. This identity holds for
every permutation quotient map, including those obtained
by inversion.

Under the involution $b\mapsto-b-11\pmod t$ from
\eqref{eq:quotient-inversion}, we may take $b=0,2,3$
when $t=6$, and $b=0,1,2,6$ when $t=8$.
For $t=8$ and $b=1,2$, the sum contains a single
nonzero term, a contradiction.

In each remaining case, exactly two indices $k,k+t$
with $0\leqslant k<t$ contribute. Hence
$\zeta=-\binom{11}{k+t}/\binom{11}{k}$.
The residues, in the order of the displayed $b$, are
\[
\begin{array}{@{}clll@{}}
\toprule
 p& b & \zeta\pmod p & \zeta^{11}\pmod p\\\midrule
 67&0,2,3&7,29,22&30,37,1\\
 89&0,6&13,71&77,55.\\
\bottomrule
\end{array}
\]
The displayed residues follow directly from the
binomial-coefficient identity above.
The only remaining case satisfying $\zeta^{11}=1$ is
$(p,b)=(67,3)$. There, $(\zeta-1)^{11}=21^{11}\equiv37\not\equiv1
\pmod{67}$, contradicting the product identity.
Thus $p=23$.

Conversely, consider $f(X)=X^{12}+2X=X(X^{11}+2)$ over
$\mathbb F_{23}$, for which $\gcd(11,22)=11$.
On nonzero squares $f(x)=3x$, whereas on nonsquares $f(x)=x$.
Since $3=7^2$ in $\mathbb F_{23}$, multiplication by $3$
permutes the nonzero squares; the nonsquares are fixed.
Together with $f(0)=0$, this proves that $f$ permutes
$\mathbb F_{23}$.
\end{proof}

\section{Characteristic classification for \texorpdfstring{$g\leqslant10$}{g <= 10}}\label{sec:low}
The factorizations of $W_g$ for $1\leqslant g\leqslant10$
are tabulated by Michon~\cite{MichonWendt}; for even $g$,
see also Fee and Granville~\cite{FeeGranville}.
For $W_g\ne0$, the following entries are exactly the sets
$\{p\text{ prime}:p\mid W_g,\ p\nmid g\}$.
Proposition~\ref{prop:wendt} makes these \emph{candidate}
characteristics, not yet the sets $\mathcal P_g$:
\[
 \begin{array}{@{}cl@{}}
\toprule
 g&\text{prime divisors}\\\midrule
 1&\varnothing\\
 2&\{3\}\\
 3&\{2,7\}\\
 4&\{3,5\}\\
 5&\{11,31\}
 \\
\bottomrule
\end{array}
 \qquad\qquad
 \begin{array}{@{}cl@{}}
\toprule
 g&\text{prime divisors}\\\midrule
 6&\text{---}\\
 7&\{2,29,127\}\\
 8&\{3,5,17\}\\
 9&\{2,7,19,37,73\}\\
 10&\{3,11,31\}
 \\
\bottomrule
\end{array}
\]
Here $W_1=1$. The dash at $g=6$ records $W_6=0$;
that case is treated separately below.

For $g\ne6$, only $(5,31)$, $(7,127)$, and $(9,73)$
remain to be excluded; all other candidates are realized
below. Since $W_6=0$, we treat $g=6$ separately.

\begin{lemma}[Three exclusions]\label{lem:three}
Let $\mathbb F_q$ be a finite field of characteristic $p$,
and let $r,s\geqslant1$ and $a\in\mathbb F_q^\times$.
If $g=\gcd(s,q-1)$ and
$(g,p)\in\{(5,31),(7,127),(9,73)\}$, then
$X^r(X^s+a)$ does not permute $\mathbb F_q$.
\end{lemma}
\begin{proof}
Suppose, to the contrary, that $f(X)=X^r(X^s+a)$
permutes $\mathbb F_q$. Put $t=(q-1)/g$ and
$\zeta=(-a)^t$. For each pair, the gcd of $Z^g-1$ and $(Z-1)^g-1$
over $\mathbb F_p$ is $Z-2$; see
Appendix~\ref{app:exclusions}. Since $g$ is odd,
$z^g=1$ implies $(z^{-1}-1)^g=-(z-1)^g$.
Thus the gcd with $(Z-1)^g+1$ is $Z-2^{-1}$.
By~\eqref{eq:prod}, this forces
$\zeta\in\{2,2^{-1}\}\subset\mathbb F_p$.
Equation~\eqref{eq:quotient-inversion} exchanges
$\zeta$ and $\zeta^{-1}$ without changing whether the
quotient map permutes $\mu_t$. Hence we may take
$\zeta=2$. Put $\tau=-g^{-1}\pmod p$ and
$Q_{g,z}(X):=P_g(\tau,X,z)$.
The modular computations give
\[
 \begin{aligned}
 \gcd_{\mathbb F_{31}[X]}(Q_{5,2},X^{31}-X)&=1,\\
 \gcd_{\mathbb F_{127}[X]}(Q_{7,2},X^{127}-X)&=1.
 \end{aligned}
\]
Appendix~\ref{app:exclusions} gives the polynomials $Q_{g,2}$
and explicit inverse polynomials for these two checks.
The third gcd is
\[
 \gcd_{\mathbb F_{73}[X]}(Q_{9,2},X^{73}-X)=X-64.
\]
Thus, for $\zeta=2$, only $U_1\equiv64\pmod{73}$
remains. Since $\tau=8$, the two values of $c$ yield
$U_2\equiv47,55\pmod{73}$, and neither gives a zero:
\begin{equation}\label{eq:73values}
 P_{18}(8,47,2)=50,\qquad
 P_{18}(8,55,2)=6\pmod{73}.
\end{equation}
Appendix~\ref{app:exclusions} records these modular
values. Inversion excludes $\zeta=2^{-1}$ as well.
These necessary congruences hold for every extension degree;
hence all three exclusions apply over arbitrary finite extensions.
\end{proof}
\begin{proposition}\label{prop:g6}
$\mathcal{P}_6=\{5,7,13,19,31\}$.
\end{proposition}
\begin{proof}
By Proposition~\ref{prop:wendt}, $p\nmid6$, so $p\ne2,3$.
Let $B(X)=(X-1)^6-1$ and choose $\omega$ with
$\omega^2+\omega+1=0$.
Since $(\omega-1)^2=-3\omega$, we have
$B(\omega)=B(\omega^2)=-28$, whereas $B(1)=-1$
and $B(-1)=63$. From
$(X^6-1)/\Phi_6(X)=(X-1)(X+1)(X^2+X+1)$,
multiplicativity of resultants gives
\[
 E_6=(-1)\cdot63\cdot(-28)^2=-2^4 3^2 7^3.
\]
Lemma~\ref{lem:residual} gives $\Phi_6(\zeta)=0$ when $p\ne7$.
For $K=6,12,18$, let $N_K(X)$ be the integer polynomial
with relatively prime coefficients and positive leading
coefficient, obtained from the rational polynomial
\[
 \operatorname{Res}_Y\bigl(\Phi_6(Y),P_K(-1/6,X,Y)\bigr)
\]
by clearing denominators and removing content.
For $j=2,3$ and $0\leqslant c<j$, let $C_{j,c}(X)$ be the
integer polynomial with relatively prime coefficients and positive
leading coefficient proportional to
\[
 N_{6j}\!\left(jX+\frac{j-1-c}{6}\right),
\]
and put $R_{j,c}:=|\operatorname{Res}_X(N_6,C_{j,c})|$.
The recurrence shows that the coefficient of $X^{K-1}$
in $P_K(T,X,Y)$ is $(Y-1)^{K-1}$.
Since $\operatorname{Res}_Y(\Phi_6,Y-1)=\Phi_6(1)=1$,
the rational resultant defining $N_K$ is monic in $X$
over $\mathbb Z[1/6]$. For $p>3$, clearing denominators
and taking primitive parts changes it only by a $p$-unit.
The affine substitutions $X\mapsto jX+(j-1-c)/6$,
$j=2,3$, have invertible slope modulo $p$, so the
corresponding vanishing conditions are preserved.

The reflected polynomial identities, with the indicated
integer normalizations, give
\[
 R_{2,1}=3^{10}R_{2,0},\qquad R_{3,2}=R_{3,0}.
\]
Only \emph{three} resultants are computed directly:
$R_{2,0}$, $R_{3,0}$, and $R_{3,1}$. The two remaining
cases follow from the reflection identities.
Their exact gcd, reproduced in Appendix~\ref{app:g6}, is
\begin{equation}\label{eq:gcd6}
 \begin{aligned}
 &\gcd(R_{2,0},R_{3,0}R_{3,1})\\[-2pt]
 &\quad=2^{220}3^{210}5^{50}7^{12}13^8 19^{13}31^5.
 \end{aligned}
\end{equation}
For $p>19$, we have $t=(q-1)/6>3$, so the inverse-power-sum
identities with indices $6,12,18$ all apply. The two possibilities for $j=2$ imply
$p\mid R_{2,0}$, and the three for $j=3$ imply
$p\mid R_{3,0}R_{3,1}$. By~\eqref{eq:gcd6}, only
$p=31$ remains.
The remaining primes $p=11,17$ are excluded by
\begin{equation}\label{eq:11-17}
 \gcd_{\mathbb{F}_p[X]}(N_6(X),X^p-X)=1
 \qquad(p=11,17),
\end{equation}
with inverse certificates in Appendix~\ref{app:g6}.
The witnesses below realize all five characteristics
$5,7,13,19,31$.
\end{proof}

Seven pairs in Theorem~\ref{thm:spectrum} have a
uniform construction.
For $q>2$, $g=q-1$, and $a\ne0,-1$, we have
\[
 X(X^{q-1}+a)=(1+a)X\qquad\text{on }\mathbb{F}_q.
\]
It yields the pairs $(g,p)=(2,3),(3,2),(4,5),
(6,7),(7,2),(8,3),(10,11)$ by taking $q=g+1$.
Table~\ref{tab:witness} gives the remaining $16$ witnesses.
For each extension field, let $\theta$ be the residue
class of $X$ modulo the irreducible polynomial listed
in the table.
For $a=\sum_{i=0}^{m-1}a_i\theta^i$ with $0\leqslant a_i<p$,
the table records $a$ as the integer $\sum_{i=0}^{m-1}a_ip^i$.
Thus, the entries in the last column are base-$p$ encodings
of field elements, not necessarily their images under the
canonical map $\mathbb Z\to\mathbb F_q$.

\begin{table}[H]
\centering\small
\caption{Witnesses not covered by the construction $g=q-1$.}
\label{tab:witness}
\begin{tabular}{rrrrlrrr}
\toprule $g$&$p$&$q$&$t$&modulus in $\mathbb{F}_p[X]$&$r$&$s$&$a$\\\midrule
3&7&7&2&$X$&1&3&3\\
4&3&9&2&$X^2+1$&1&4&3\\
5&11&11&2&$X$&1&5&2\\
6&5&25&4&$X^2-2$&1&6&7\\
6&13&13&2&$X$&1&6&2\\
6&19&19&3&$X$&1&6&5\\
6&31&31&5&$X$&5&6&5\\
7&29&29&4&$X$&1&7&4\\
8&5&25&3&$X^2+3$&1&8&1\\
8&17&17&2&$X$&1&8&3\\
9&2&64&7&$X^6+X+1$&1&9&8\\
9&7&343&38&$X^3+5$&1&9&171\\
9&19&19&2&$X$&1&9&5\\
9&37&37&4&$X$&4&9&9\\
10&3&81&8&$X^4+X+2$&1&10&21\\
10&31&31&3&$X$&1&10&11\\
\bottomrule
\end{tabular}
\end{table}
For each row, the accompanying script checks the
irreducibility of the field modulus (if its degree exceeds $1$),
the identities $\gcd(s,q-1)=g$ and $\gcd(r,g)=1$,
and the bijectivity of
$y\mapsto y^r(y^{s/g}+a)^g$ on $\mu_t$.
It evaluates this quotient map at exactly $91$ subgroup
elements across the $16$ rows; see Appendix~\ref{app:witnesses}.
The candidate lists from $W_g$, Lemma~\ref{lem:three},
Proposition~\ref{prop:g6}, and the $23$ witnesses now
prove Theorem~\ref{thm:spectrum}.

\FloatBarrier
\begin{corollary}[Classification for seven parameter pairs]
\label{cor:fan-seven}
Let
\[
 (g,p)\in\{(2,3),(3,2),(4,5),(6,7),(7,2),(8,3),(10,11)\}.
\]
Suppose $q=p^m$, $r,s\geqslant1$, $a\in\mathbb F_q^\times$, and
$\gcd(s,q-1)=g$. Write $q_0:=g+1=p^d$. Then $d\mid m$; put
\[
 n:=m/d,\qquad t:=(q-1)/g,\qquad u:=s/g.
\]
For $n=1$, the binomial $X^r(X^s+a)$ permutes $\mathbb F_q$
if and only if $\gcd(r,g)=1$ and $a\ne-1$.
For $n\geqslant2$, it permutes $\mathbb F_q$ if and only if
\[
 \gcd(r,g)=1,\qquad (-a)^t\ne1,\qquad
 r\frac{q_0^h-1}{q_0-1}\equiv u\pmod t
\]
for some $1\leqslant h<n$ with $\gcd(h,n)=1$.
\end{corollary}
\begin{proof}
The divisibility $p^d-1=g\mid p^m-1$ forces $d\mid m$.
Moreover, $q-1=(q_0-1)t$ and $s=(q_0-1)u$,
so $\gcd(u,t)=1$. When $n\geqslant2$, apply
\cite[Corollary~6.1]{Fan26} over $\mathbb F_{q_0}$,
with extension degree $n$ and exponent index $u$.
When $n=1$, the polynomial agrees with $(1+a)X^r$
on $\mathbb F_{q_0}$, giving the stated criterion.
\end{proof}

\appendix
\section{Reproducible algebraic certificates}\label{app:cert}
We record the finite certificates for
Theorem~\ref{thm:spectrum}. Appendix~\ref{app:exclusions}
handles the three exclusions, Appendix~\ref{app:g6}
records the $g=6$ calculations, and
Appendix~\ref{app:witnesses} verifies the permutation witnesses.
Theorem~\ref{thm:finite} is independent of these computations.

\subsection{Three exclusions}\label{app:exclusions}
The three modular gcds used in Lemma~\ref{lem:three} are
\[
\begin{array}{@{}cl@{}}
\toprule
 (g,p)&\gcd_{\mathbb F_p[Z]}(Z^g-1,(Z-1)^g-1)\\\midrule
 (5,31)&Z-2\\
 (7,127)&Z-2\\
 (9,73)&Z-2\\
\bottomrule
\end{array}
\]
The negative-sign cases follow by inversion, as proved
in Lemma~\ref{lem:three}, and need no further computation.
For $(g,p)=(5,31),(7,127)$, let $I_{g,p}$ be an inverse
of $X^p-X$ modulo $Q_{g,2}$ in $\mathbb F_p[X]$.
The following coefficient vectors are listed in ascending
powers of $X$:
\[
 \begin{array}{@{}ll@{}}
\toprule
 \text{object}&\text{computed coefficient vector}\\\midrule
 Q_{5,2}&[18,19,17,4,1]\\
 I_{5,31}&[22,0,13,3]\\
 Q_{7,2}&[50,111,47,31,95,17,1]\\
 I_{7,127}&[71,24,69,90,17,122]
 \\
\bottomrule
\end{array}
\]
The modular identities
$(X^{31}-X)I_{5,31}\equiv1\pmod{Q_{5,2}}$ and
$(X^{127}-X)I_{7,127}\equiv1\pmod{Q_{7,2}}$
prove the two gcds are $1$. For $(9,73)$, the three
computations needed for the orientation $\zeta=2$ are
\[
\begin{array}{@{}lr@{}}
\toprule
\text{object}&\text{result in }\mathbb F_{73}\\\midrule
\gcd(Q_{9,2},X^{73}-X)&X-64\\
P_{18}(8,47,2)&50\\
P_{18}(8,55,2)&6\\
\bottomrule
\end{array}
\]
The last two values are precisely the nonzero residues
used in~\eqref{eq:73values}.
The inverse orientation follows from~\eqref{eq:quotient-inversion}.

\subsection{The case \texorpdfstring{$g=6$}{g=6}}\label{app:g6}
The polynomials $N_6,N_{12},N_{18}$ are defined in the
proof of Proposition~\ref{prop:g6}. The script constructs
them exactly and computes only the three integer
resultants $R_{2,0}$, $R_{3,0}$, and $R_{3,1}$.
Their gcd is~\eqref{eq:gcd6}; the two other cases follow
from the reflection identities, without further computations.

The coefficients of $N_6(X)$, in descending order, are
\begin{quote}\footnotesize\ttfamily
60466176, -1763596800, 22515252480,\\
-165288211200, 770586032448, -2376500847840,\\
4892764590720, -6614063618400, 5595135918876,\\
-2665119096360, 542601329725.
\end{quote}
For the two congruences in~\eqref{eq:11-17}, the
inverses modulo $N_6$, listed in ascending order, are
\[
 \begin{array}{@{}cl@{}}
\toprule
 p&I_p\\\midrule
 11&[2,6,0,10,5,3,10,1,4,1]\\
 17&[0,4,11,0,1,8,8,7,10,15].
 \\
\bottomrule
\end{array}
\]
Thus $(X^p-X)I_p\equiv1\pmod{N_6}$ in $\mathbb F_p[X]$
for $p=11,17$.

\subsection{Verification of permutation witnesses}\label{app:witnesses}
For the $16$ rows of Table~\ref{tab:witness}, the script
constructs the indicated finite fields and checks the
irreducibility of each extension-field modulus. It also
checks $\gcd(s,q-1)=g$ and $\gcd(r,g)=1$.
Writing $u=s/g$ and $t=(q-1)/g$, it then forms
$\mu_t\subset\mathbb F_q^\times$ and verifies the equality
\[
 \{y^r(y^u+a)^g:y\in\mu_t\}=\mu_t.
\]
This is precisely the quotient criterion in
Lemma~\ref{lem:quotient}. The $16$ subgroups contain a
total of $91$ elements. The other seven characteristic
pairs are realized by the construction preceding Table~\ref{tab:witness}.

The accompanying single-file script
\texttt{PBpWendt\_arXiv1-verify.py} uses Python~3 and SymPy
for exact arithmetic in three verification modules: the modular exclusions of Appendix~\ref{app:exclusions},
the $g=6$ certificates of Appendix~\ref{app:g6},
and the permutation witnesses above. It needs no external
data files or project scripts. The cited factorizations of
$W_g$, the proved value of $E_6$, and
Corollary~\ref{cor:primefield-g11} do not require computer
verification. Nor does the assertion $73\nmid D_9$ in
Section~\ref{sec:primefield}, which follows from the
explicit gcd and product given there.

\providecommand{\bysame}{\leavevmode\hbox to3em{\hrulefill}\thinspace}
\providecommand{\MR}{\relax\ifhmode\unskip\space\fi MR }
\providecommand{\MRhref}[2]{%
  \href{http://www.ams.org/mathscinet-getitem?mr=#1}{#2}
}
\providecommand{\href}[2]{#2}


\begin{thebibliography}{ABK15}

\bibitem[ABK14]{ABK14}
Mohamed Ayad, Kacem Belghaba, and Omar Kihel, \emph{On permutation binomials
  over finite fields}, Bull. Aust. Math. Soc. \textbf{89} (2014), no.~1,
  112--124,
  \href{https://doi.org/10.1017/S0004972713000208}{\nolinkurl{doi:10.1017/S0004972713000208}}.

\bibitem[ABK15]{ABK15}
\bysame, \emph{On permutation binomials}, Rocky Mountain J. Math. \textbf{45}
  (2015), no.~2, 389--399,
  \href{https://doi.org/10.1216/RMJ-2015-45-2-389}{\nolinkurl{doi:10.1216/RMJ-2015-45-2-389}}.

\bibitem[AW06]{AkbaryWang}
Amir Akbary and Qiang Wang, \emph{A generalized {Lucas} sequence and
  permutation binomials}, Proc. Amer. Math. Soc. \textbf{134} (2006), no.~1,
  15--22,
  \href{https://doi.org/10.1090/S0002-9939-05-08220-1}{\nolinkurl{doi:10.1090/S0002-9939-05-08220-1}}.

\bibitem[Fan26]{Fan26}
Xiang Fan, \emph{A complete classification of permutation binomials of the form
  {$X^r(X^{q-1}+a)$} over finite fields}, 2026,
  \href{https://arxiv.org/abs/2609.20354v1}{\nolinkurl{arXiv:2609.20354v1}}.

\bibitem[FG91]{FeeGranville}
Greg Fee and Andrew Granville, \emph{The prime factors of {Wendt}'s binomial
  circulant determinant}, Math. Comp. \textbf{57} (1991), no.~196, 839--848,
  \href{https://doi.org/10.1090/S0025-5718-1991-1094948-8}{\nolinkurl{doi:10.1090/S0025-5718-1991-1094948-8}}.

\bibitem[HPL23]{HouPL}
Xiang-dong Hou and Vincenzo Pallozzi~Lavorante, \emph{New results on
  permutation binomials of finite fields}, Finite Fields Appl. \textbf{88}
  (2023), 102179,
  \href{https://doi.org/10.1016/j.ffa.2023.102179}{\nolinkurl{doi:10.1016/j.ffa.2023.102179}}.

\bibitem[Leh35]{Lehmer35}
Emma Lehmer, \emph{On a resultant connected with {Fermat}'s last theorem},
  Bull. Amer. Math. Soc. \textbf{41} (1935), no.~12, 864--867,
  \href{https://doi.org/10.1090/S0002-9904-1935-06210-X}{\nolinkurl{doi:10.1090/S0002-9904-1935-06210-X}}.

\bibitem[Mic]{MichonWendt}
G{\'e}rard~P. Michon, \emph{Factorization of {Wendt}'s determinant},
  Numericana, online table,
  \href{https://www.numericana.com/data/wendt.htm}{\nolinkurl{www.numericana.com/data/wendt.htm}}.
\bibitem[MZ07]{MZ07}
Ariane~M. Masuda and Michael~E. Zieve, \emph{Nonexistence of permutation
  binomials of certain shapes}, Electron. J. Combin. \textbf{14} (2007), no.~1,
  N12, \href{https://doi.org/10.37236/1013}{\nolinkurl{doi:10.37236/1013}}.

\bibitem[MZ09]{MZ}
\bysame, \emph{Permutation binomials over finite fields}, Trans. Amer. Math.
  Soc. \textbf{361} (2009), no.~8, 4169--4180,
  \href{https://doi.org/10.1090/S0002-9947-09-04578-4}{\nolinkurl{doi:10.1090/S0002-9947-09-04578-4}}.

\bibitem[PL98]{ParkLee}
Young~Ho Park and June~Bok Lee, \emph{Permutation polynomials with exponents in
  an arithmetic progression}, Bull. Aust. Math. Soc. \textbf{57} (1998), no.~2,
  243--252,
  \href{https://doi.org/10.1017/S0004972700031622}{\nolinkurl{doi:10.1017/S0004972700031622}}.

\bibitem[Sim02]{Simalarides02}
Anastasios Simalarides, \emph{Upper bounds for the prime divisors of {Wendt}'s
  determinant}, Math. Comp. \textbf{71} (2002), no.~237, 415--427,
  \href{https://doi.org/10.1090/S0025-5718-00-01292-8}{\nolinkurl{doi:10.1090/S0025-5718-00-01292-8}}.

\bibitem[Wan87]{Wan87}
Daqing Wan, \emph{Permutation polynomials over finite fields}, Acta Math.
  Sinica (N.S.) \textbf{3} (1987), no.~1, 1--5,
  \href{https://doi.org/10.1007/BF02564938}{\nolinkurl{doi:10.1007/BF02564938}}.

\bibitem[Wen94]{Wendt}
Ernst Wendt, \emph{{Arithmetische Studien \"uber den ``letzten'' Fermatschen
  Satz, welcher aussagt, dass die Gleichung $a^n=b^n+c^n$ f\"ur $n>2$ in ganzen
  Zahlen nicht aufl\"osbar ist}}, J. Reine Angew. Math. \textbf{113} (1894),
  335--347,
  \href{https://doi.org/10.1515/crll.1894.113.335}{\nolinkurl{doi:10.1515/crll.1894.113.335}}.

\bibitem[WL91]{WanLidl}
Daqing Wan and Rudolf Lidl, \emph{Permutation polynomials of the form {$x^r
  f(x^{(q-1)/d})$} and their group structure}, Monatsh. Math. \textbf{112}
  (1991), no.~2, 149--163,
  \href{https://doi.org/10.1007/BF01525801}{\nolinkurl{doi:10.1007/BF01525801}}.

\bibitem[Zie08]{Zieve}
Michael~E. Zieve, \emph{Some families of permutation polynomials over finite
  fields}, Int. J. Number Theory \textbf{4} (2008), no.~5, 851--857,
  \href{https://doi.org/10.1142/S1793042108001717}{\nolinkurl{doi:10.1142/S1793042108001717}}.


\end{thebibliography}
\end{document}